\documentclass[12pt,a4paper,oneside,onecolumn]{article}
\title{Classical Poisson algebra of a vector bundle : Lie-algebraic characterization}
\author{Lecomte P.B.A. and Zihindula Mushengezi E.}
\usepackage{etex}
\usepackage[]{latexsym}
\usepackage[T1]{fontenc}
\usepackage[latin1]{inputenc}
\usepackage[english]{babel}
\usepackage{graphics}
\usepackage[]{amssymb}
\usepackage{multicol}
\usepackage{amsmath}
\usepackage{graphicx}
\usepackage{amsfonts}
\usepackage[all]{xy}
\usepackage{pifont}
\usepackage{microtype}
\usepackage{answers}
\usepackage{tabularx}
\usepackage{float}
\usepackage{color}

\usepackage{listings}
\definecolor{darkWhite}{rgb}{0.94,0.94,0.94}
\usepackage{palatino}

\newcount\exenum     \newdimen\numindent
\newcount\qqnum     \newdimen\qqmarge
\newcount\sqnum     \newdimen\sqmarge

\outer\def\bye{%
\vskip 0pt\@endmulticol\@endgroup              
\ifanswer \let\next=\exobye@                   
\else     \let\next=\@exobye                   
\fi\next}

\usepackage[Glenn]{fncychap}

\usepackage{pstricks}
\usepackage{pgf,tikz}
\usepackage{pstricks-add}
\usepackage{pst-plot}
\usetikzlibrary{arrows}
\usepackage{pgfplots}
\usepackage{pgfkeys}

\definecolor{qqqqqq}{rgb}{0,0,0}

\definecolor{xdxdff}{rgb}{0.49,0.49,1}
\definecolor{qqwuqq}{rgb}{0,0.39,0}
\definecolor{qqqqff}{rgb}{0,0,1}
\definecolor{ttttff}{rgb}{0.2,0.2,1}
\definecolor{uququq}{rgb}{0.25,0.25,0.25}

\renewcommand{\epsilon}{ \varepsilon}

\newcommand{\euro}{\texteuro{}}

\newcommand{\pre}{{\bf Proof.\ }}

\newcommand{\cl }{\mathcal }

\newtheorem{theo}{Theorem}[section]
\newtheorem{prop}[theo]{Proposition}

\newcounter{exercice}

\definecolor{fbase}{rgb}{0.8,0.8,1}
\definecolor{fgris}{gray}{0.6}
\definecolor{frouge}{HTML}{DC143C}
\definecolor{fvert}{rgb}{0.6,1,0.6}
\definecolor{fbleu}{rgb}{0.4,0.4,1}
\definecolor{fjaune}{HTML}{DCDC14}
{\endMakeFramed}
\makeatletter
\DeclareRobustCommand\sfrac[1]{\@ifnextchar/{\@sfrac{#1}}%
                                            {\@sfrac{#1}/}}
\def\@sfrac#1/#2{\leavevmode\kern.1em\raise.5ex
         \hbox{$\m@th{\fontsize\sf@size\z@\selectfont#1}$}
         \kern-.1em/\kern-.15em\lower.55ex
          \hbox{$\m@th{\fontsize\sf@size\z@\selectfont#2}$}}

\DeclareRobustCommand{\Efrac}[2]{{\displaystyle\begingroup
\raise2ex\hbox{$\m@th{#1}$}\endgroup\@@over \lower1ex
\hbox{$\m@th{#2}$}}}
\makeatother

\numberwithin{equation}{section}

\newtheorem{Exc}{Exercice}
\Newassociation{correction}{Soln}{mycor}
\Newassociation{indication}{Indi}{myind}

\def\exo#1{\futurelet\testchar\MaybeOptArgmyexoo}
\def\MaybeOptArgmyexoo{\ifx[\testchar \let\next\OptArgmyexoo
                        \else \let\next\NoOptArgmyexoo \fi \next}
\def\OptArgmyexoo[#1]{\begin{Exc}[#1]\normalfont}
\def\NoOptArgmyexoo{\begin{Exc}\normalfont}

\newcommand{\finexo}{\end{Exc}}
\newcommand{\flag}[1]{}

\tikzset{
xmin/.store in=\xmin, xmin/.default=-3, xmin=-3,
xmax/.store in=\xmax, xmax/.default=3, xmax=3,
ymin/.store in=\ymin, ymin/.default=-3, ymin=-3,
ymax/.store in=\ymax, ymax/.default=3, ymax=3,
}

\newcommand{\entete}[1]

\begin{document}
 \Opensolutionfile{mycor}[ficcorex]
 \Opensolutionfile{myind}[ficind]
 \entete{\'Enoncés}


\maketitle

\begin{abstract}
We prove that the Lie algebra $\cl S(\cl P(E,M))$ of symbols of linear operators acting on smooth sections of a vector bundle $E\to M,$ characterizes it. To obtain this, we assume that $\cl S(\cl P(E,M))$ is seen as ${\rm C}^\infty(M)-$module and that the vector bundle is of rank $n>1.$ 

We improve this result for the Lie algebra $\cl S^1(\cl P(E,M))$ of symbols of first-order linear operators. We obtain a Lie algebraic characterization of vector bundles with $\cl S^1(\cl P(E,M))$ without the hypothesis of being seen as a ${\rm C}^\infty(M)-$module.  
\end{abstract}

\section{Introduction}

We know that the space of all differential operators acting on the sections of a vector bundle is a quantum Poisson algebra as is the case for differential operators acting on smooth functions. See in \cite{LecZih} for instance.

But for differential operators acting on sections, as we will see in this article, the study of symbols gives rise to two interpretations:
\begin{enumerate}
\item[$\bullet$] The usual principal symbol linked to the order of derivation
\item[$\bullet$] The symbol related to the filtration of the structure of quantum Poisson algebra.
\end{enumerate}
This particularity is part of the richness of the Poisson algebra structure of the symbols of these differential operators and gives rise to remarkable developments.

Let $ E\to M$ be a vector bundle of rank $n.$ Denoting by $\Gamma(E)$ the space of smooth sections of $E,$ let us consider the Lie algebra of linear operators of this vector bundle defined by  
\[
   \cl D(E,M) =\bigcup_{k\geq 0} \cl D^k(E,M),
  \]
where
  \[
    \cl A(E,M):=\cl D^0(E,M)=\{T\in End(\Gamma(E)):[T,\gamma_u]=0,\,\forall u\in{\rm C}^\infty(M)\}
  \]
and for any integer $k\geq1,$
  \[  
    \cl D^k(E,M)= \{T\in End(\Gamma(E))| \forall u\in{\rm C}^\infty(M): [T,\gamma_u]\in    
    \cl D^{k-1}(E,M)\}.
  \]
Provided with the previous filtration, the Lie algebra $\cl D(E,M)$ is not a quantum Poisson algebra. But we have the result below which we will use later.
\begin{prop}\label{D quasi dist}
The Lie algebra $\cl D(E,M)$ is quasi-distinguishing; that is, the relations
 \begin{enumerate}
  \item $([T,Z(\cl A)]=0)\Rightarrow T\in \cl A(E,M)$ 
  \item $\{T\in\cl D(E,M):[T, Z(\cl A)]\subset \cl D^k(E,M)\}=\cl D^{k+1}(E,M),$  $\forall k\in\mathbb{N}$
 \end{enumerate}
 are both satisfied, where $Z(\cl A)=\{\gamma_u, u\in{\rm C}^\infty(M))\}$ is the center of the associative algebra $\cl A(E,M).$
\end{prop}

In the rest of this article, we have assumed that the rank $ n $ of the vector bundle $ E $ exceeds 1. In fact, Graboswski and Poncin established in \cite{GraPon4} that for $ n = 1, $ the Lie algebras $\cl D (E, M) $ and $\cl D (M),$ the last algebra being that of all differential operators acting on ${\rm C}^\infty(M),$ are isomorphic. 
It is clear that this does not allow us to consider any characterization of a vector bundle with such a Lie algebra.\\
In \cite{LecZih} we have obtained Lie-algebraic characterization with $\cl D(E,M),$ for vector bundles of rank  greater than 1. We have also provided this Lie algebra with appropriated filtration in order to make it a quantum Poisson algebra. \\
With this new filtration, $\cl D(E,M)$ becomes $\cl P(E,M)$ and we can now define the space of symbols of differential operators in $\cl P(E,M),$ as in \cite{GraPon1,GraPon2,GraPon4}, for instance.\\

This is what we discuss in the following lines.

\section{The classical Poisson algebra $\cl S(\cl P(E,M))$}

Let $ E \to M $ be a vector bundle of rank $ n> 1. $ 
The space $ \Gamma (E) $ of smooth sections being a $ {\rm C}^\infty (M) - $ module, let us recall that
  \[
  \gamma_u:\Gamma(E)\to\Gamma(E):s\mapsto us, \forall u\in{\rm C}^\infty(M)
  \]
is an endomorphism of the space $ \Gamma (E). $ \\  
We pose 
\[
  \cl P(E,M)=\bigcup_{k\geq 0}\cl P^k(E,M),
\]
with, by definition,
\[
  \cl P^0(E,M)=\{\gamma_u: u\in{\rm C}^\infty(M)\}
\] 
and, for any integer $k,$
\[
 \cl P^{k+1}(E,M)=\{T\in End(\Gamma(E))|\forall u\in{\rm C}^\infty(M): [T,\gamma_u]\in\cl P^k(E,M)\}.
 \]

The following results are taken from \cite{LecZih}.
\begin{prop}\label{pem quantum}
For any integer $j,k\in\mathbb{N},$ we have
\begin{enumerate}
  \item $\cl P^k(E,M)\subset \cl P^{k+1}(E,M)$ and $\cl P^{j}(E,M)\cdot\cl P^{k}(E,M)\subset\cl P^{j+k}(E,M)$
  \item $[\cl P^{j}(E,M),\cl P^{k}(E,M)]\subset\cl P^{j+k-1}(E,M)$
\end{enumerate}
\end{prop}
The above Proposition \ref{pem quantum} means the Lie algebra $\cl P(E,M)$ is a quantum Poisson one.
\begin{prop}
The quantum Poisson algebra $ \cl P (E, M) $ satisfies the following properties.
\begin{enumerate}
\item $[T,\gamma_u]=0,\forall u\in{\rm C}^\infty(M)\Rightarrow T\in gl(E)\subset \cl P^1(E,M)$
\item $\{T\in\cl P(E,M)| [T,\cl P^0(E,M)]\subset \cl P^{k}(E,M)\}=\cl P^{k+1}(E,M),\quad$ $\forall k\in\mathbb{N}$
\end{enumerate}
\end{prop}
The first property above is the  non-singularity of $\cl P(E,M)$. This last algebra is also sympletic; this meaning that its center $Z(\cl P(E,M))$ contains only constants, i.e, multiplication by elements of $\mathbb{R}.$\\
We also have the following surprising equality
\[
\cl D(E,M)=\cl P(E,M).
\]

In what follows, we study the classical limit of the quantum Poisson algebra $ \cl P (E, M).$ 

Let us specify that the  classical limit mentioned above is defined by
  \[
     \cl S(\cl P(E,M))=\bigoplus_{i\in\mathbb{Z}}\cl S^i(\cl P(E,M));
  \]  
with $\cl S^i(\cl P(E,M))=\cl P^i(E,M)/\cl P^{i-1}(E,M).$
We obtain a classical Poisson algebra whose operations are given in the following.\\
Let us recall that for any $T\in\cl P^i(E,M),$ $ord(T)=i,$ if $T\notin \cl P^{i-1}(E,M).$\\
For $i\geq ord(T),$ the $i$-degree symbol of $ T $  is defined by
\[
\sigma_i(T)=\left\{
\begin{array}{l}
 0 \mbox{ if } i>ord(T)\\
 T+\cl P^{i-1}(E,M)\mbox{ if } i=ord(T)
\end{array}
\right.
\]
The symbol related to the quantum Poisson structure of $ \cl P (E, M)$ is for its part given by
  \[
  \sigma_{pson} : \cl P(E,M)\to  \cl S(\cl P(E,M)): T\mapsto \sigma_{ord}(T).
  \]
For $P\in \cl S^i(\cl P(E,M))$ and $Q\in\cl S^j(\cl P(E,M))$ such that $P=\sigma_i(T)$ and $Q=\sigma_j(D),$ we set, by definition,
 \[
  P.Q=\sigma_{i+j}(T\circ D) \mbox{ and } \{P,Q\}=\sigma_{i+j-1}([T,D])\cdot
 \]
 
\section{Particular case of  the Lie sub-algebra $gl(E)\subset\cl P^1(E,M)$}\label{cas part sle}
 
By virtue of the calculations made in the previous section, we have, by definition, 
 \[
 \sigma_{pson}(\gamma_u)=\gamma_u+\{0\},\quad \forall u\in{\rm C}^\infty(M)
 \]
and we will simply denote $\qquad\sigma_{pson}(\gamma_u)=\gamma_u.$\\

Likewise, for $A\in gl(E)\setminus\cl P^0(E,M),$ we have
\[
 \sigma_{pson}(A)= A'+\cl P^0(E,M),
\] 
with $A'=A-\frac{tr(A)}{n}id,$\quad $tr(A)$ being the trace of $A.$ 
Therefore, for any $A,B\in gl(E),$ we have 
\begin{eqnarray*}
 \sigma_{pson}([A,B]) & = & [A,B]+ \cl P^0(E,M)\\
                      & = & [A',B']+\cl P^0(E,M).
\end{eqnarray*} 
And for the product, if $A,B\notin\cl P^0(E,M),$ we have 
  \[
  \sigma_{pson}(A)\cdot\sigma_{pson}(B)=0,
  \]
but if $\gamma_u\in\cl P^0(E,M),$ we then have   
  \[
  \sigma_{pson}(\gamma_u)\cdot\sigma_{pson}(A)=\gamma_u \circ A'+\cl P^0(E,M).
  \]
We therefore have the following identification of Lie algebras
\[
\sigma_{pson}(gl(E))\cong sl(E)\oplus {\rm C}^\infty(M)id,
\]
where the multiplication is commutative and defined by
 \[
 (A+\gamma_u)\cdot(B+\gamma_v)=\gamma_v\circ A +\gamma_u\circ B +\gamma_{uv};
 \] 
and the bracket being given by the following relation
\[
  \{A+\gamma_u , B+\gamma_v \}=[A,B].
\]

\section{General case}

Let us state the following result which gives the local expression of the elements of $ \cl P (E, M) $ in a trivialization of $ E. $

Sometimes, for convenience of writing, we simply denote $ \gamma_u $ by $ u \in {\rm C}^\infty (M).$
\begin{prop}\label{critere pkde }
The elements of $\cl P^k(E,M),$ $k\geq1,$ are characterized by the fact that they are written locally, in a trivialization domain $ U \subset M, $ in the form
\begin{equation}\label{(*v)}
  \sum_{|\alpha|<k}T_\alpha\partial^\alpha+\sum_{|\beta|=k}u_\beta\partial^\beta
\end{equation}
where $T_\alpha\in{\rm C}^\infty (U,gl(n,\mathbb{R}))$ and $u_\beta\in{\rm C}^\infty (U).$
\end{prop}
\pre
The proof is done by induction on $ k $. Let $ T \in \cl P^1(E, M) $ and suppose that in a domain of trivialization $ U $ of $ E $ we have
  \[
    T=A+\sum_i^mT_i\partial_i,\quad A,T_i\in{\rm C}^\infty (U,gl(n,\mathbb{R}))
  \]    
The relation $[T,\gamma_u]\in\cl P^0(E,M)$ gives 
 \[
  \sum_{i=1}^m T_i\circ\partial_i(u)\in{\rm C}^\infty(U),\quad\forall u\in{\rm C}^\infty(M).
 \]
Therefore, for any $i\in[1,m]\cap\mathbb{N},$ we do have  $T_i\in{\rm C}^\infty(U).$ 
 
Assume by induction that the result is true for any element of $\cl P^r(E,M)$ with $r<k.$ 
Let $T\in\cl P^k(E,M)$ such that we locally have
 \[
 T=\sum_{|\alpha|<k }T_\alpha\partial^\alpha+\sum_{|\beta|=k}T_\beta\partial^\beta  .
\] 
We have, by applying the induction hypothesis to $ [T, \gamma_u], $ that its highest order term of derivation is of the form
 \[
 \sum_{|\lambda|=k}u_\lambda\partial^\lambda.
 \]
Now in $ [T_\beta \partial^\beta, \gamma_u], $ the highest order terms of derivation are of the form $\partial_iu T_{\beta}\partial^{\beta_i},$ with $\beta_i=\beta-e_i,$ where $(e_i)_{1\leq i\leq m}
$ designates the canonical basis of the $\mathbb{R}-$vector space $\mathbb{R}^m.$ \\ 
We therefore have, for any $u,$ $\partial_{i}uT_{\beta}\in{\rm C}^\infty(U),$ and consequently $T_{\beta}\in{\rm C}^\infty(U).$

Conversely if $ T \in\cl D (E, M) $ is written locally in the form (\ref{(*v)}), we obtain that $T\in\cl P^k(E,M),$ for all $k\geq 1,$ by doing a recurrence on $k.$

Thus, in a trivialization domain $ U \subset M, $ we obtain the part of order strictly equal to $ k $ of the local expression of $ T.$\\
It indeed has the following form
\begin{equation}\label{(**)}
 \sum_{|\alpha|=k-1}A_{\alpha}\partial^\alpha+\sum_{|\beta|=k}u_{\beta}\partial^\beta, 
\end{equation}
where $A_\alpha\in {\rm C}^\infty(U,sl(n,\mathbb{R}))$  and $u_\beta\in{\rm C}^\infty(U).$\\
Particularly, for any $D\in\cl P^l(E,M), $ having   
\[
 \sum_{|\lambda|=l-1}B_{\lambda}\partial^\lambda+\sum_{|\mu|=l}v_\mu\partial^\mu  \]
as $l-$order terms,  those of $k+l-$order of $T\circ D$ (but also those of $D\circ T$) are therefore given by 
\begin{equation}\label{produit loc}
   \sum_{|\alpha|=k-1}\sum_{|\mu|=l}A_{\alpha}v_\mu\partial^\alpha\partial^\mu
   +
  \sum_{|\lambda|=l-1}\sum_{|\beta|=k} B_\lambda u_\beta\partial^\beta\partial^\lambda
   +
    \sum_{|\beta|=k}\sum_{|\mu|=l}u_\beta v_\mu\partial^\beta\partial^\mu
\end{equation}
\newpage
For the Lie bracket $[T,D],$ the $k+l-1-$order terms are  
\begin{equation}\label{crochet loc}
  \sum_{|\alpha|=k-1}\sum_{|\lambda|=l-1}[A_\alpha,B_\lambda]\partial^\alpha\partial^\lambda
\end{equation}
\[  
+
   \sum_{|\alpha|=k-1}\sum_{|\mu|=l}\sum_{1\leq i\leq m}\left(A_\alpha\partial_iv_\mu\partial^{\alpha_i}\partial^\mu
    - \partial_iA_\alpha v_\mu\partial^\alpha\partial^{\mu_i} \right)
\]
\[    
+ 
    \sum_{|\beta|=k}\sum_{|\lambda|=l-1}\sum_{1\leq i\leq m}\left(u_\beta\partial_iB_\mu\partial^\lambda\partial^{\beta_i}
-\partial_iu_\beta B_\lambda\partial^\beta\partial^{\lambda_i}\right)
\]
\[
+\sum_{|\beta|=k}\sum_{|\lambda|=l-1}\sum_{1\leq i\leq m}\left(u_\beta\partial_iv_\mu\partial^{\beta_i}\partial^\mu
-   
       v_\mu\partial_ju_\beta\partial^{\mu_j}\partial^\beta\right)
\] 

We notice that the local decomposition (\ref{(**)})  given above is not intrinsic.
In fact, if in the sum on the right we recognize the principal symbol of the differential operator in the usual sense, the sum on the left, for its part, does not resist a change of coordinates and is therefore not globally defined.

In the following lines we build a global decomposition allowing to find a global meaning to the expression given in (\ref{(**)})  previously.\\

Let now $T\in\cl S^{k-1}(M)\otimes sl(E)$ and assume given a unit partition $(U_i,\rho_i)$ of $M$ whose domains $U_i$ are the trivialisation one of $E.$ In any $U_i,$ if $T$ is expressed in the form
 \[
   T=\sum_{|\alpha|=k-1}A_{\alpha,i}\xi^\alpha\cdot
 \] 
We then set
 \[
  \overline{T}_i=\sum_{|\alpha|=k-1}A_{\alpha,i}\partial^\alpha\in\cl D^{k-1}(E|_{U_i},U_i)
 \]
with $A_{\alpha,i}\in  {\rm  C}^\infty(U_i,sl(n,\mathbb{R})).$ 
The differential operator
  \[
    \overline{T}=\sum_i \rho_i\overline{T}_i\in \cl D^{k-1}(E,M)\subset\cl P^k(E,M),
  \]
associated with the partition of the unit chosen at the start is then such that
   \[
    \sigma_{pson}(\overline{T})=\sigma_{ppal}(\overline{T})=T, \footnote{Note that $\sigma_{ppal}$ is the usual principal symbol.}
   \]
but it is obviously not the only one to verify this relation. \\
Nevertheless, we have the following statement.
\begin{prop}
The space $\cl S^k(\cl P(E,M))=\cl P^k(E,M)/\cl P^{k-1}(E,M)$ of symbols in the sense " quantum Poisson algebra " of differential operators in $\cl P^k(E,M)$ is determined by the following  exact short sequence of $\mathbb{R}-$vector spaces 
   \[
    0\longrightarrow\cl S^{k-1}(M)\otimes sl(E)\stackrel{\theta}{\longrightarrow}\cl P^k(E,M)/\cl P^{k-1}(E,M)\stackrel{\delta}{\longrightarrow}\cl S^k(M)\longrightarrow 0,
   \] 
with $\theta: T\mapsto \overline{T}+\cl P^{k-1}(E,M)$ and 
  \[
  \delta: D+\cl P^{k-1}(E,M)\mapsto \left\{
                                 \begin{array}{l}
                                 0 \mbox{ if } D\in\cl D^{k-1}(E,M)\\
                                 \sigma_{ppal}(D) \mbox{ if not. }
                                 \end{array}
                                \right. 
  \] 
\end{prop}
\pre
The application $\theta$ is well-defined. Indeed, the differential operators $D_1,D_2\in\cl D^{k-1}(E,M)\subset\cl P^k(E,M)$ are such that $\sigma_{pson}(D_1)=\sigma_{pson}(D_2),$ we then do have $D_1-D_2\in\cl P^{k-1}(E,M);$ which means that the image of $ T $ does not depend on the choice of the operator $\overline{T}$ such that $\sigma_{pson}(\overline{T})=T.$ \\ 
Also, $ \theta $ is obviously a linear map and it is injective. Indeed, let $T\in  S^{k-1}(M)\otimes sl(E)$ such that $\theta(T)=0.$ 
We then have $\overline{T}\in\cl P^{k-1}(E,M).$ 
But by construction $\overline{T}\notin\cl D^{k-2}(E,M).$ We thus have
    \[
      \sigma_{ppal}(\overline{T})\in S^{k-1}(M)id
    \] 
We deduce, since $\sigma_{pson}(\overline{T})=\sigma_{ppal}(\overline{T})=T,$ that $T=0.$  \\
 The map $ \delta $ being linear and directly surjective, we show to finish that
\[
ker(\delta)=Im(\theta).
\] 
The inclusion $ker(\delta)\supset Im(\theta)$ is obvious. 
Let's prove the other sense of that inclusion. If $D+\cl P^{k-1}(E,M)\in ker(\delta),$ then, by the definition of $\delta,$ we have
 \[
   D\in \cl D^{k-1}(E,M)\cap\cl P^k(E,M).
 \]	
Consequently, on the one hand,
 \[
  \sigma_{ppal}(D)\in \cl S^{k-1}(M)\otimes gl(E),
 \] 
because $D\in\cl D^{k-1}(E,M)$ and on the other, like $ D \in \cl P^k (E, M), $ we have instead,
\[
   \sigma_{ppal}(D)\in \cl S^{k-1}(M)\otimes sl(E).
\] 
The inclusion sought is a direct result of this.\hfill $\blacksquare$\\
\newpage
Therefore, seen as $\mathbb{R}-$vector spaces, we have the following decomposition 
  \[
    \cl S^k(\cl P(E,M))=Pol^{k-1}(T^*M,sl(E))\oplus Pol^k(T^*M,\mathbb{R})
  \]
for any integer   $k\in\mathbb{N}.$ \\ 
 
The question that arises is whether the following exact sequence of Lie algebras (but also of  associative algebras), whose exactness comes from that given in the previous statement and from the operations previously performed in (\ref{produit loc}) and in (\ref{crochet loc}) is split. 
\begin{equation}\label{(v)}
 \xymatrix{
 \quad 0\ar[r] & \cl S(M)\otimes sl(E)\ar[r] & \cl S(\cl P(E,M))\ar[r] & \cl S(M)\ar[r]& 0 
         }
\end{equation}
Note that the splitting of this sequence leads in particular to that of the following exact sequence of Lie algebras
\[
 \xymatrix{
 \quad 0\ar[r] & sl(E)\ar[r] & \cl S^1( \cl P(E,M))\ar[r] & Vect(M)\ar[r]& 0 
         }
\]

To answer it, we make use of the following result where a split sequence of Lie algebras is given.
\begin{prop}\label{suite scindee}
Let $E\to M$ be vector bundle of rank $n.$  
With respect to a connection on $ E, $ the following short exact sequence of Lie algebras is split
\[
\xymatrix{       
  0\ar[r] & \cl P^0(E,M)\stackrel{i}{\longrightarrow} \cl P^1(E,M)\stackrel{\sigma_{pson}}{\longrightarrow} \cl P^1(E,M)/\cl P^0(E,M)\ar[r] & 0 
         }
\]
where $i$ is the canonical injection and $\sigma_{pson}$ the map previously defined.
\end{prop}

\pre
Via a covariant derivation $ \nabla $ of $ E, $ we have the identification of $ \mathbb {R}-$vector spaces 
  \[
    \cl P^1(E,M)\cong Vect(M)\oplus gl(E).
  \]
Indeed, for $T\in\cl P^1(E,M)\subset\cl D^1(E,M),$ we have 
\[
  \sigma_{ppal}(T)=X\in Vect(M),
\] 
with the identification $\cl S^1(M)id\simeq Vect(M).$\\

Therefore, since $\nabla_X\in\cl P^1(E,M),$ the difference $T-\nabla_X$ is an endomorphisms field. 

We denote by $\lambda: \cl P^1(E,M)\to Vect(M)\oplus gl(E)$ the linear bijection thereby defined. If we have 
 $ T=\nabla_X+A$ and $D=\nabla_Y+B,$ we then get in $\cl P^1(E,M)$
  \[
   [T,D]=\nabla_{[X,Y]}+R^\nabla(X,Y)+ \nabla_XB-\nabla_YA+[A,B]\cdot
  \] 
This is equivalent to
 \begin{equation}\label{(*)}
  \quad [(X,A),(Y,B)]=([X,Y],R^\nabla(X,Y)+ \nabla_XB-\nabla_YA+[A,B])
 \end{equation}
in the space $Vect(M)\oplus gl(E).$   

Moreover, consider the following short exact sequence
\[
\xymatrix{       
  0\ar[r] & sl(E)\ar[r] & \cl P^1(E,M)/\cl P^0(E,M)\ar[r] & Vect(M)\ar[r] &  0 
         }
\]
corresponding to the particular case $ k = 1 $ of that characterizing the space of symbols, in the sense "Quantum Poisson", of the differential operators of order $ k$, that we gave at the beginning of this section.

 Thus, as in the general case at the beginning of this section, we have the linear map
   \[
  \delta: T+\cl P^{0}(E)\mapsto \left\{
                                 \begin{array}{l}
                                 0 \mbox{ if } T\in gl(E)\\
                                 \sigma_{ppal}(T) \mbox{ if not }
                                 \end{array}
                                \right. 
  \] 
which is surjective and the injection $\theta: A\in sl(E)\mapsto A+\cl P^0(E,M).$

Seen as $\mathbb{R}-$vector spaces, we therefore have the following identification
 \[
    \cl P^1(E,M)/\cl P^0(E,M)\cong Vect(M)\oplus sl(E). 
 \]
Consider the following commutative diagram
\[
  \xymatrix{(X,A)\in Vect(M)\oplus sl(E) \ar[d]^{\lambda^{-1}}\ar[r]^\mu & \cl P^1(E,M)/\cl P^0(E,M)\\  
  \nabla_X+A\in\cl P^1(E,M)\ar[ur]_{\sigma_{pson}} &
           }
\]
The linear map $ \mu = \sigma_{pson} \circ \lambda^{- 1} $ is injective because $ \mu (X, A) = 0 $ induces $ \nabla_X + A \in \cl P^0 (E); $ and we deduce that
\[
X=0 \quad\mbox{ and } \,A\in sl(E)\cap\cl P^0(E,M).
\]   
Likewise, $ \mu $ is surjective.\\
The bracket in $\cl P^1(E,M)/\cl P^0(E,M)$ is given by
\[
   [[T],[D]]  = \nabla_{[X,Y]}+R^\nabla(X,Y)+ \nabla_XB-\nabla_YA+[A,B]+\cl P^0(E,M)
 \]
with $[T]=\nabla_X+A+\cl P^0(E,M)$ and $[D]=\nabla_Y+B+\cl P^0(E,M).$
 
Therefore, the corresponding operation in $ Vect (M) \oplus sl (E), $ obtained by structure transport via $ \mu, $ is not necessarily a Lie bracket since the term $ R^\nabla ( X, Y) $ is not always of zero trace. \\

To remedy this, suppose that $ \nabla $ is associated with a connection form of a reduction of the frame principal bundle $ L^1(E) $ of $ E $ to the Lie subgroup $ O (n) $ of $ GL (n, \mathbb {R}). $ For such a derivation, $ R^\nabla $ has values in $ sl (E). $\\
We then have an isomorphism of Lie algebras
\[
\mu : Vect(M)\oplus sl(E)\to \cl P^1(E,M)/\cl P^0(E,M),
\]
the space $ Vect (M) \oplus sl (E) $ being provided with the following bracket
 \begin{equation}\label{(***)}
   [(X,A),(Y,A)]=([X,Y],R^\nabla(X,Y)+ \nabla_XB-\nabla_YA+[A,B]).
 \end{equation}  
Given the relations  (\ref{(*)})  and  (\ref{(***)}), we conclude that the canonical injection
\[
 \beta: Vect(M)\oplus sl(E)\to Vect(M)\oplus gl(E)
\]
is a homomorphism of Lie algebras.
Consequently,
\[
  \lambda^{-1}\circ \beta\circ \mu^{-1}:\cl P^1(E,M)/\cl P^0(E,M)\to \cl P^1(E,M)
\]
is a homomorphism of Lie algebras allowing to identify the Lie algebra $ \cl P^1 (E, M) / \cl P^0 (E, M) $ to a Lie subalgebra of $ \cl P^1 (E, M), $ for the structures specified in the previous lines, and we can see that it is a section of $ \sigma_{pson}. $ \\
We have just shown that the short exact sequence of the statement is split. \\
\hspace*{4cm}\hfill $\blacksquare$

Note that if $ \nabla $ is a covariant derivation of $ E $ associated with the reduction in question in the previous proof, then for $ T = \nabla_X + A, $ the following decomposition
\[
T=(\nabla_X+A-\frac{1}{n}tr(A)\,id)+\frac{1}{n}tr(A)\,id\,\cdot
\]
only depends on the reduction and not on the choice of connection.\\ 
Indeed, if relatively to another covariant derivation $ \nabla', $ associated with the same reduction, we consider an analogous decomposition of $ T, $ then we have
 \begin{eqnarray*}
   T=\nabla'_X+A'& = & \nabla_X+(A'+(\nabla'_X-\nabla_X))\\
                 & = & \nabla_X+S+A' 
 \end{eqnarray*}
with $S=\nabla_X-\nabla'_X,$ and thus the trace of $S$ is null. We conclude that $A'=A-S$. \\
 Therefore, $tr(A)=tr(A')$ and we have
  \[
   \nabla'_X+A'-\frac{1}{n}tr(A')\,id=\nabla_X+S+A-S-\frac{1}{n}tr(A)\,id.
  \] 

Now consider the following diagram 
  \[
  \xymatrix{           &    & 0\ar[d]           &   0\ar[d] \\
                       &    &  gl(E)\ar[d]      &  sl(E)\ar[d]       &    \\                   
\quad 0\ar[r] & \cl P^0(E,M)\ar[r] & \cl P^1(E,M)\ar[r]\ar[d] & \cl P^1(E,M)/\cl P^0(E,M)\ar[r]\ar[d]& 0\\ 
                       &    &  Vect(M)\ar[d]    & Vect(M)\ar[d]& \\
                       &    &  0                & 0 
         }
  \]

We have established that the horizontal sequence is split. And we deduce that the split of the vertical sequence on the right would lead to that of the vertical sequence located on the left. Now according to \cite{Lec suit ex}, the splitting of this sequence essentially requires the naturality of the vector bundle $ E. $\\
Note that a vector bundle $E\to M$ is necessarily natural if the base $ M $ is simply connected.\\
 
Thus, the answer to the question of whether the exact sequence of Lie algebras given in $ (\ref{(v)}) $ is always split is negative.

\section{Lie-algebraic characterization of vector bundles}

Let us begin by stating results of Lie-algebraic characterization of vector bundles taken from \cite{Lec3}. 

\begin{theo}\label{ende carct e}
Let $ E \mapsto M $ and $ F \mapsto M $ be two vector bundles of respective ranks $n,n'>1$ with $H^1(M,\mathbb{Z}/2)=0$. The Lie algebras $ gl (E) $ and $ gl (F) $ (resp. $ sl (E) $ and $ sl (F) $) are isomorphic if and only if the vector bundles $ E $ and $ F $ are isomorphic.
\end{theo}

We use the above Theorem \ref{ende carct e} to obtain the following result.
\begin{theo}
Let $E\to M, F\to M$ be two vector bundles of respective ranks $n,n'>1$ with $H^1(M,\mathbb{Z}/2)=0.$
The Lie algebras $ \cl S (\cl P (E, M)) $ and $ \cl S (\cl P (F, M)), $ seen as $ {\rm C}^\infty (M) -$modules, are isomorphic if, and only if, the vector bundles $ E $ and $ F $ are.
\end{theo}
\pre
Let $\Phi : \cl S (\cl P (E, M))\to \cl S (\cl P (F,M)) $ be an isomorphism of Lie algebras.
Note that $\Phi$ preserves $\cl A={\rm C}^\infty(M)$ the basis of these classical Poisson algebras.
Indeed, for any $u\in {\rm C}^\infty(M),$ we have
\[
  \Phi(\gamma_u)=u\Phi(\gamma_1)
\]
and the conclusion comes from the fact that the quantum Poisson algebras considered are symplectic.\\
Let $B\in gl(E).$ We then have
\[
  \{B,\cl A \}=0.
\]
Therefore $\{\Phi(B),\cl A\}=0$ and this implies $\Phi(B)\in gl(E),$ since $\cl D(E,M)$ is quasi-distinguishing, according to the Proposition \ref{D quasi dist}. 
\hfill $\blacksquare$\\

For the Lie subalgebras $ \cl S^1 (\cl P (E, M)) $ and $ \cl S^1 (\cl P (F, N)), $ this result may improve. To prove this, we use the short exact sequence presented in the previous Proposition \ref{suite scindee}.

\begin{theo}
Let $E\to M, F\to M$ be two vector bundles of respective ranks $n,n'>1$ with $H^1(M,\mathbb{Z}/2)=0.$
The Lie algebras $ \cl S (\cl P^1(E, M)) $ and $ \cl S (\cl P^1(F, M)), $  are isomorphic if, and only if, the vector bundles $ E $ and $ F $ are. 
\end{theo}

\pre
We observe, starting from the decomposition
\[
\cl S^1(\cl P(E,M))=sl(E)\oplus Vect(M)
\] 
obtained previously via a connection on $ E.$ \\ 

For all $ T = (\nabla_X, A),$ let $ B \in sl (E) $ such that there exists $ r \in \mathbb {N} $ verifying 
\[
 (ad(T))^r(B)=0.
\]
We then have 
 \[
   \nabla_X(\nabla_X\cdots(\nabla_X(B)))=0,
 \]
where $\nabla_X$ is $r$ times applied.
 In a trivialization of $E$ of domain $U\subset M,$ considering $B$ whose local expression is of the form $(\alpha_{ij})=(\delta_{12}u),$ $u\in {\rm C}^\infty(U),$ i.e. having all its terms null except that of the position $(1,2),$ we can choose $u$ so that we necessarily have $X=0.$ 
We deduce that
  \[
   Nil(\cl S^1(\cl P(E))\subset sl(E).
  \] 
Moreover, we know that if $A\in sl(E)$ is such that $A^p=0,$ we then have $A\in Nil(\cl S^1(P(E,M)).$ 
Indeed, observe that
\begin{eqnarray*}
 (ad A)^k(\nabla_X+B) & = & (ad A)^k(\nabla_X)+(ad A)^k(B) \\
                      & = & -(ad A)^{k-1}(X\cdot A)+ (ad A)^k(B)
\end{eqnarray*} 
where, by virtue of the particular case studied in the previous section \ref{cas part sle}, for all $C,D\in sl(E),$ we have
\begin{eqnarray*}
  (ad\, C)(D) & = & \{C,D\} =  [C,D] ,       
\end{eqnarray*}
which allows us to conclude, since $ (ad \, C)^k (D) $ is then a sum of the terms of the form $a_k C^k\circ D\circ C^{k-1}, a_k\in\mathbb{R}$.

In the following lines, the goal is to establish that the linear envelope of the nilpotent endomorphism fields is the entire space $ sl (E) $. Let $ A \in sl (E). $ Over a trivialization domain $ U \subset M, $ we can therefore write
  \[
    A|_U=\Sigma_{i}N_i^U  
  \]
with $N_i^U \in sl(E|_U),  (1\leq i\leq n^2-1)$, which are nilpotent endomorphism fields since $sl(n, \mathbb {R}) $ admits a basis formed of nilpotent matrices.
Consider now a Palais cover of $M,$ 
\[
  \cl O=\cl O_1\cup\cdots\cup\cl O_r, r\in \mathbb{N},
\]
locally finite, the elements $ U _{\alpha,j} $ of each $ \cl O_j $ being trivialization domains of $ E $ 2 by 2 disjoint, and  $(\cl \rho_{\alpha,j}),$ a partition of the unit, locally finite and subordinate to this cover. 

We therefore have 
\begin{eqnarray*}
  \rho_{\alpha,j}A & = & \sum_{i=1}^{n^2-1} N_{i,\alpha,j}
\end{eqnarray*} 
with $N_{i,\alpha,j}\in sl(E|_{U_{\alpha,j}})$ nilpotent and compactly supported in $U_{\alpha,j}$. \\
Let now pose $\cup\cl O_j=\cup_{\alpha}U_{\alpha,j}=U_j$ and consider $N_{ji}$ defined by
\[
     N_{ji}(x)=  \left\{ \begin{array}{l}
           N_{i,\alpha,j}(x) \mbox{ if } x\in U_{j}\\
           0 \quad\quad\quad \mbox{ if not }
               \end{array} 
               \right.
\] 
We then obtain the smoothness of $N_{ji}$. Indeed, for any $x\notin U_j,$ consider an open neighborhood $V\ni x$ of compact adherence. We know that $ V $ is only encountered by a finite number of supports of $N_{\alpha,j},$ whose reunion is the compact that we agree to denote by $K.$ 

Therefore $V\setminus K$ is an open neighborhood of $ x $ in which $ N_{ji} $ is identically zero. And we have thus established our assertion since for $ x \in U_j $ the smoothness of  $N_{ji}$ is obvious. 
We conclude that 
\[
  \sum_{\alpha} \rho_{\alpha,j} A=\sum_{i=1}^{n^2-1} N_{ji},
\] 
and as a result, 
 \[
   A=\sum_{i=1}^{n^2-1}\sum_{j=1}^rN_{ji}.
 \]
We have therefore just shown that
\[
\left.  \right\rangle  Nil(\cl S^1(\cl P(E,M))\left\langle= sl(E), 
  \right.
\]
where the usual notation $\left.  \right\rangle  H\left\langle\right.$ designates the linear envelope of the subset $H$ of a vector space.

We deduce that for any isomorphism 
 \[
 \Phi: \cl S(\cl P^1(E,M))\to  \cl S(\cl P^1(F,N))
 \]
of Lie algebras is necessarily such that
\[
  \Phi(sl(E))= sl(F).
\]
Hence, by virtue of the previous  Theorem \ref{ende carct e}, we have the desired result. \\
\hspace*{1cm}\hfill $\blacksquare$


\end{document}